\documentclass[twoside,11pt]{article}

\usepackage{blindtext}

\usepackage{jmlr2e}
\usepackage{hyperref}
\usepackage{url}
\usepackage{wrapfig}
\usepackage{graphicx}
\usepackage{subcaption}
\usepackage{xspace}
\usepackage{float} 
\usepackage{array}
\usepackage{makecell}
\usepackage{natbib}

\usepackage{amsmath,amsfonts}
\usepackage{algorithmic}
\usepackage{algorithm}
\usepackage{amssymb}
\usepackage{tabularx}

\usepackage{lastpage}

\begin{document}

\title{Uncertain but Useful: Leveraging CNN Training Variability into Data Augmentation}

\author{
  \name Inés Gonzalez-Pepe
  \email inesgp99@gmail.com \\
  \addr Department of Computer Science\\
  Concordia University\\
  Montreal, QC  H3G 1M8, Canada \\
  \AND
  \name Vinuyan Sivakolunthu
\email vinuyansivakolunthu@gmail.com \\
  \addr Department of Computer Science\\
  Concordia University\\
  Montreal, QC  H3G 1M8, Canada \\
  \AND
  \name Yohan Chatelain
  \email yohan.chatelain@camh.com \\
  \addr Krembil Centre for Neuroinformatics \\
  Centre for Addiction and Mental Health \\
  Toronto, ON M5T 0S8, Canada \\
  \AND
  \name Tristan Glatard
  \email tristan.glatard@camh.ca \\
  \addr Krembil Centre for Neuroinformatics \\
  Centre for Addiction and Mental Health \\
  Toronto, ON M5T 0S8, Canada \\
  }

\maketitle

\begin{abstract}
  Deep learning (DL) has transformed neuroimaging by delivering state-of-the-art performance with reduced computation times. Yet, the numerical uncertainty inherent to DL training remains largely underexplored despite its potential to significantly impact the reliability of model outcomes. We show that training the FastSurfer segmentation model introduces substantial numerical uncertainty that exceeds its non-DL counterpart (FreeSurfer 7.3.2) in cortical regions, potentially impacting downstream clinical results. We also characterize this training-time uncertainty using random seed perturbations and demonstrate that seed-induced variability is structurally comparable to numerical variability. We then show that seed variability can be leveraged as a data augmentation technique through ensembling to improve downstream brain age regression performance. These findings position numerical uncertainty during DL training as a substantive factor in neuroimaging reliability, with measurable consequences for downstream tasks, and demonstrate that it can simultaneously be harnessed as a data augmentation technique.
\end{abstract}

\begin{keywords}
  Numerical uncertainty, Convolutional neural networks, Deep ensembles, Data augmentation
\end{keywords}

\section{Introduction}
Deep learning (DL) has transformed neuroimaging, yet the numerical stability of the training phase remains unexplored. Unlike inference, training involves iterative, high-dimensional optimization whose trajectory is sensitive to small numerical perturbations that accumulate across iterations, potentially steering models toward different solutions even on identical data.
While numerical uncertainty is well-documented in traditional pipelines~\citep{glatard2015reproducibility,alizadeh2025numerical,salari2021accurate,kiar2021data}, understanding its behavior during DL training is critical as these neural networks increasingly replace classical toolkits.  This work addresses this critical gap and presents two major findings. 

First, using Monte Carlo Arithmetic (MCA)~\citep{parker1997monte}, a well-established numerical technique, we present the first numerical uncertainty quantification of a neuroimaging DL model during training, demonstrating that FastSurfer's~\citep{henschel2022fastsurfervinn}~training-time variability meets or exceeds that of its traditional counterpart, FreeSurfer~\citep{fischl2002whole}, particularly in cortical structures. This finding is directly tied to recent work by~\citep{chatelain2026practical}, which measured and formalized how comparable numerical noise propagates through statistical estimators to distort clinical findings.


Second, we establish that random seed perturbations serve as a highly practical and validated proxy for numerical variability. While MCA provides a theoretically grounded method to track floating-point precision errors, its exhaustive execution at runtime is computationally prohibitive. We demonstrate empirically that random seed variations generate regional variability patterns that are structurally comparable to MCA, providing a low-cost framework to measure training uncertainty natively. Having observed that models initialized with different seeds diverge in anatomically valid ways, we reframe this training-time instability through ensembling as a data augmentation technique. 
We show that this approach performs comparably to traditional synthetic data generators in downstream brain age regression across both internal and external datasets.
All code required to reproduce this work is publicly available at \url{https://anonymous.4open.science/r/cnn_training_variability-EBB7}.

\section{Materials \& Methods}

\subsection{Model Perturbations}

\paragraph{Monte Carlo Arithmetic.}
To simulate numerical perturbations (such as hardware and software variability) during training, we apply MCA, a stochastic arithmetic method that introduces random perturbations in
floating-point operations at runtime by applying the following perturbation:
\begingroup
\[  random\_rounding(x \circ y) = round(inexact(x \circ y)) \]
\endgroup
where $x$ and $y$ are floating-point numbers, $\circ$ is an arithmetic operation, and $inexact$ is a random perturbation defined at a given virtual precision:
\begingroup
\[ inexact(x)=x+2^{e_{x}-t}\xi \]
\endgroup
where \(e_{x}\) is the exponent in the floating-point representation of \(x\),
\(t\) is the virtual precision, and \(\xi \sim U(-\frac{1}{2}, \frac{1}{2})\) is
a uniform random variable.
To measure numerical uncertainty, we applied a perturbation of 1 ulp (unit of least precision), which corresponds to a virtual precision of $t=24$ bits for single-precision and $t=53$ bits for double-precision.

\paragraph{Fuzzy PyTorch.} To apply MCA during FastSurfer training, we used Fuzzy PyTorch~\citep{fuzzy-pytorch}, a framework that integrates stochastic arithmetic directly into PyTorch's execution pipeline. Fuzzy PyTorch combines the Verificarlo compiler~\citep{denis2015verificarlo} with the PRISM backend to introduce controlled floating-point perturbations at runtime without modifying model code or architecture. Specifically, we applied stochastic rounding mode and trained FastSurfer on CPU, as MCA is currently supported only on CPU.

\paragraph{MCA in FreeSurfer.} We use FreeSurfer v7.3.2, in order to match the version used by the FastSurfer authors during the original model training. 
To measure the variability in FreeSurfer whole brain segmentation, we applied MCA to FreeSurfer using ``fuzzy libmath"~\citep{salari2021accurate}, a version of the GNU mathematical library instrumented with the Verificarlo compiler and ran FreeSurfer on CPU.

\paragraph{Random Seeds.}
To provide a variability baseline and a practical proxy for MCA, we trained each model on GPU using 15 different seeds to assess the impact of initialization and stochastic training dynamics.
In PyTorch, random seeds control pseudorandom
processes during training, including weight initialization and dropout. To
isolate model uncertainty (a.k.a. epistemic uncertainty), data shuffling was
disabled and the dataset order fixed.

\subsection{FastSurfer Training} 
The FastSurfer whole-brain segmentation CNN is composed of three 2D fully convolutional neural networks, each associated with a different 2D slice orientation. 
We retrain the FastSurfer model v2.4.0 and follow the authors' original protocol using the code available on GitHub~\citep{fastsurfer-github} (including using a base learning rate of 0.01, the "cosineWarmRestarts" learning rate  scheduler and AdamW optimizer) and train it on  FreeSurfer-generated segmentations, with all inputs preprocessed using the FreeSurfer HiRes stream with the
Desikan–Killiany–Tourville (DKT) atlas. 
We obtained access to the original
training and validation datasets (illustrated in Table~\ref{table:datasets} with data acknowledgements on Github), replacing the
Rhineland Study (restricted access) with the MICA dataset of comparable
submillimeter resolution. All datasets listed were used to train and evaluate the model, except for CoRR and HBN, which were reserved for downstream experiments. As in the original paper, all MRI data was conformed to a standard orientation and 1mm isotropic resolution prior to input into the model. No additional preprocessing was performed.
After, all segmentations were quality-controlled by visually inspecting each image in sagittal, coronal, and axial views across multiple slice locations to identify anatomical errors or imaging artifacts. 

We then compare the CNN with the
classical FreeSurfer pipeline on a subset of the Consortium for Reliability and Reproducibility (CoRR) dataset~\citep{zuo2014open} that aims to evaluate test-retest reliability and reproducibility, and has previously been used for evaluating variability during FastSurfer inference~\citep{pepe2023numerical}.

To assess the impact of uncertainty on segmentations, we
compute the minimum Sørensen-Dice score across pairs of repetitions (MCA or random seed). Given $N$ different repetitions of a given segmentation of a brain region $S$, let $S_i$ and $S_j$
represent the set of voxels assigned to $S$ the $i$-th
and $j$-th repetitions, respectively. The uncertainty of the segmentation is evaluated as:
\[
  U_S = \min_{\substack{i, j \in \{1, \dots, N\} \\ i \neq j}} \left( \frac{2\  |S_i \cap S_j|}{|S_i| + |S_j|} \right)
\]
where $|S_i \cap S_j|$ is the number of overlapping voxels classified as part
of $S$ in both repetitions, and $|S_i| + |S_j|$ are the total
number of voxels assigned to $S$ in each repetition. Taking the minimum across all pairs highlights the worst-case consistency between repetitions, making it a conservative measure of variability. 


\begin{table}[tbp]
\centering
\begin{tabular}{c >{\centering\arraybackslash}p{1.4cm} >{\centering\arraybackslash}p{1.4cm} >{\centering\arraybackslash}p{1.4cm} >{\centering\arraybackslash}p{1.4cm}} 
 \hline

 Dataset & Number of Subjects & Voxel Size (mm) & Age Range (years) & Brain Age Experiment \\ [0.5ex] 
 \hline
HCP~\citep{VanEssen2013HCP} & 129 & 0.7 & 22-35 & FastSurfer \\
ABIDE-I~\citep{di2014autism} & 20 & 1.0 & 18-64 & FastSurfer \\ 
ABIDE-II~\citep{di2017enhancing} & 25 & 0.9 & 20-31 & FastSurfer \\ 
ADNI~\citep{Mueller2005ADNI} & 63 & 1.0 & 55-93 & FastSurfer \\ 
IXI~\citep{ixi_dataset} & 43 & 1.0 & 19-87 & FastSurfer \\ 
LA5C~\citep{poldrack2016phenome} & 40 & 1.0 & 21-50 & FastSurfer \\ 
MIRIAD~\citep{Malone2012MIRIAD} & 12 & 1.0 & 55-86 & FastSurfer \\ 
OASIS1~\citep{Marcus2007OASIS1} & 55 & 1.0 & 18-90 & FastSurfer \\ 
OASIS2~\citep{Marcus2009OASIS2} & 28 & 1.0 & 60-96 & FastSurfer \\ 
MICA~\citep{royer2022open} & 50 & 0.8 & N/A & N/A \\ 
CoRR~\citep{zuo2014open} & 32 & 1.0 & N/A & N/A \\ 
HBN~\citep{alexander2017open} & 1000 & 1.0 & 5-22 & HBN \\ 
 \hline
\end{tabular}
\caption{Overview of datasets used to train and evaluate FastSurfer models and to conduct downstream data augmentation experiments. The HBN dataset served as an external test set for assessing generalization in brain age prediction.}
\label{table:datasets}
\end{table}

\subsection{FastSurfer Ensembling}

We exploit FastSurfer’s numerical variability as an augmentation mechanism. Segmentations from models trained under different perturbations remain anatomically valid, enabling numerical ensembling at a one-time training cost.

{Random seeds capture comparable epistemic uncertainty as MCA through initialization and stochastic dynamics, at a fraction of the computational cost, making them a practical augmentation mechanism.
Multiple FastSurfer models were trained with different random seeds, and the same inputs were processed through each model. An ensemble of 15 seed models therefore produced 15 distinct yet valid segmentations per subject. Region of interest (ROI) volumes (in $mm^3$) were then computed from voxel counts across segmentations and used in downstream analysis.
Unlike conventional deep ensembling~\citep{lakshminarayanan2017simple}, which aggregates predictions at inference, we leverage training variability as an anatomically realistic augmentation technique.
Brain age prediction served as the evaluation task, using ROI volumes as features for Random Forest (RF), Support Vector Machine (SVM), and Gradient Boosting (GB) regression models. 
Generalization was further evaluated on the unseen Healthy Brain Network (HBN) dataset~\citep{alexander2017open}.

To contextualize our approach, we compared numerical ensembling to synthetic data augmentation methods. We first explored variational autoencoders, which achieved moderate MAE but produced unrealistic ROI distributions, likely due to limited training data; these results were therefore excluded from quantitative comparisons. A Gaussian Copula synthesizer~\citep{sdv_synthetic_data_vault} generated more plausible anatomical distributions and was used as the primary augmentation baseline. To prevent data leakage, synthetic generators were trained exclusively on the training split, and MAE was evaluated on held-out test sets.


\section{Results}

\subsection{CNN Numerical Variability Exceeds That of Traditional Method}
\label{sec:var_results}

\begin{figure}[tbp]
    \centering
    \includegraphics[width=0.99\linewidth]{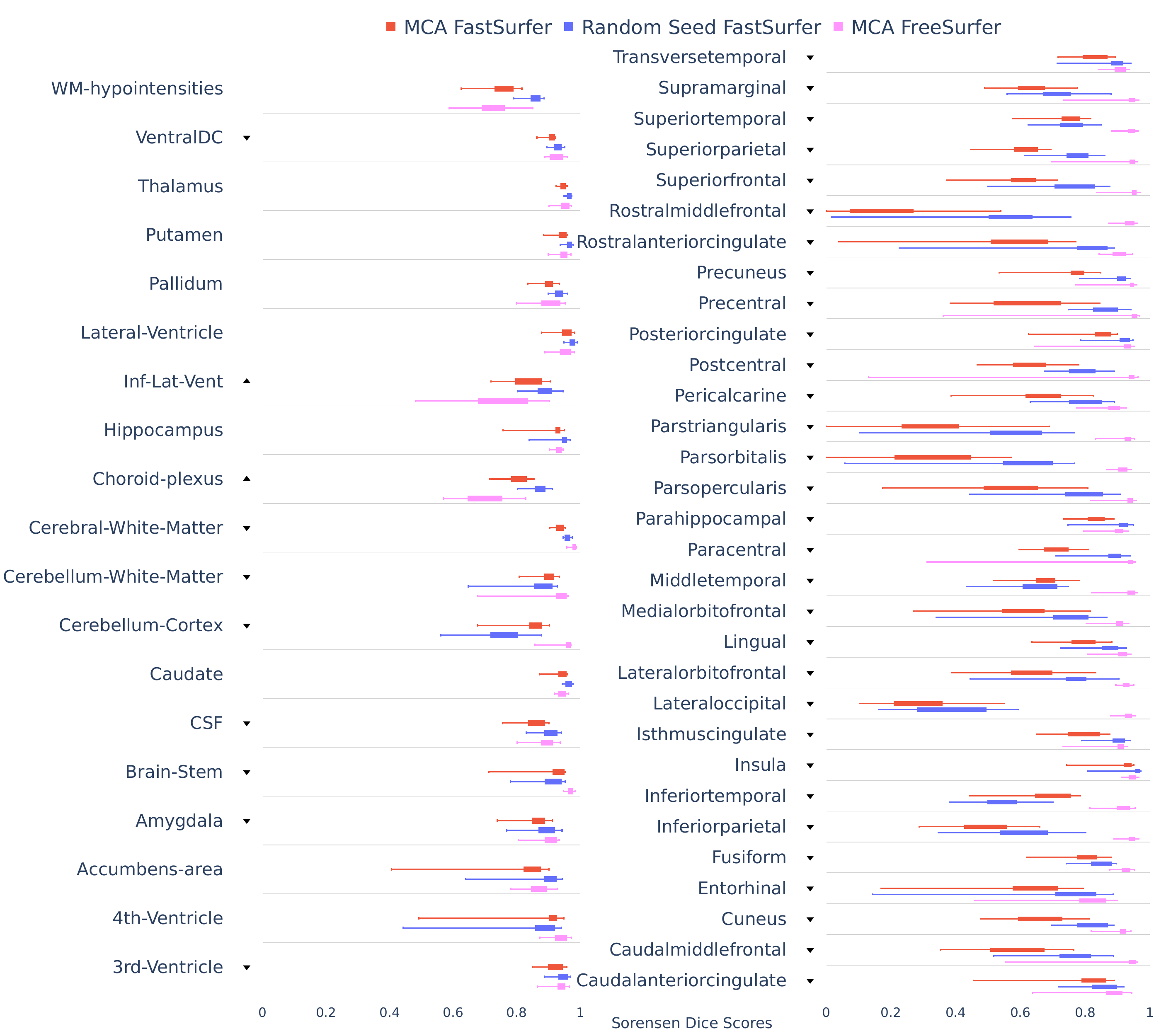}
    \caption{Minimum Sørensen-Dice variability ($U_S$) across FastSurfer and FreeSurfer ROI on CoRR data (32 subjects). $\triangle$ indicate regions where FastSurfer is less
    variable, and $\triangledown$ where it is more variable than FreeSurfer
    (one-sided t-test, p<0.05, Bonferroni corrected). Left: subcortical regions. Right: cortical regions. 
    }
    \label{fig:random_comparison_min_dice}
\end{figure}
\begin{figure}[tbp]
    \centering

    \begin{subfigure}[t]{0.33\textwidth}
        \centering
        \includegraphics[width=\linewidth]{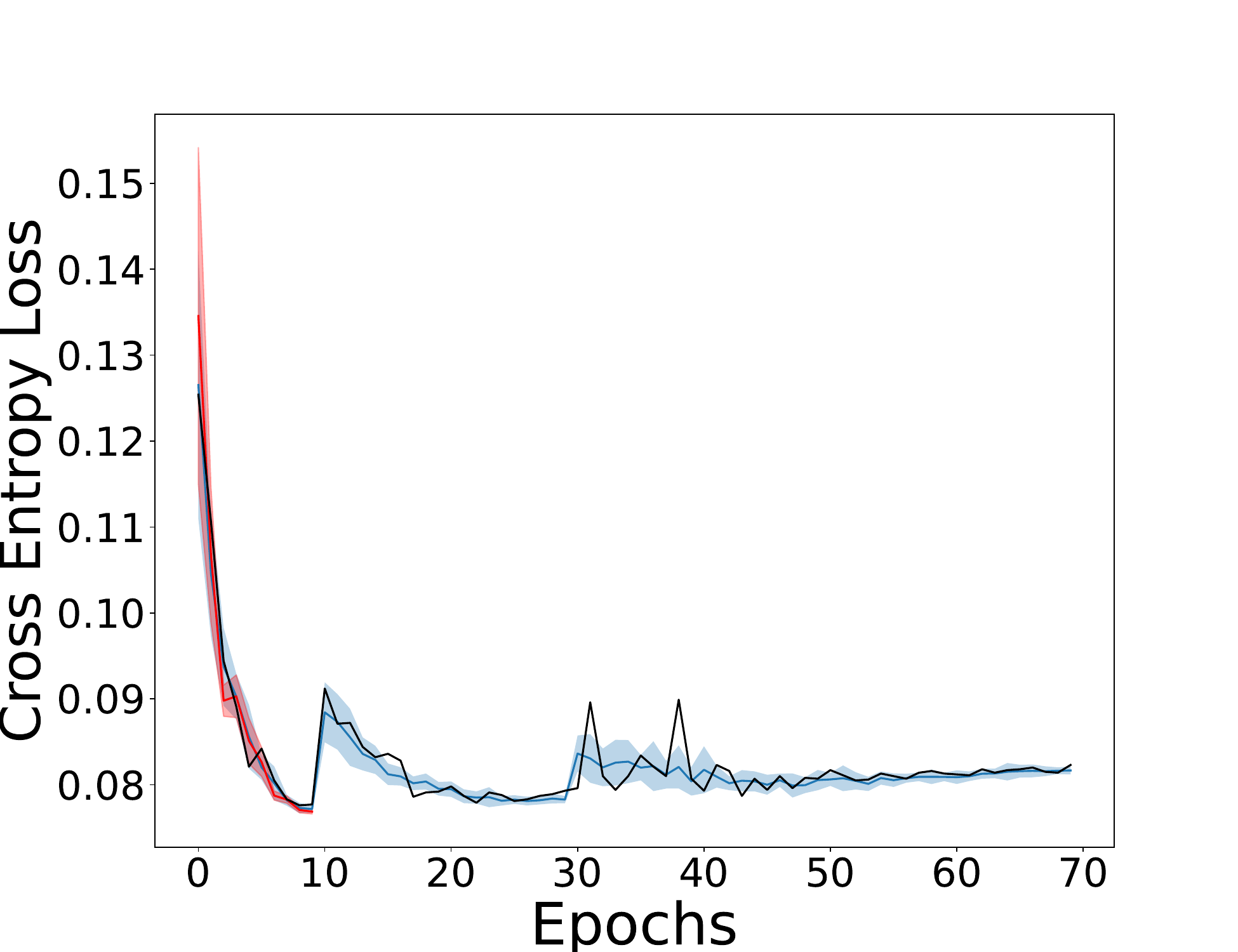}
        \caption{Cross-Entropy Loss}
        \label{fig:loss_ce}
    \end{subfigure}
    \hspace{-1em}
    \begin{subfigure}[t]{0.33\textwidth}
        \centering
        \includegraphics[width=\linewidth]{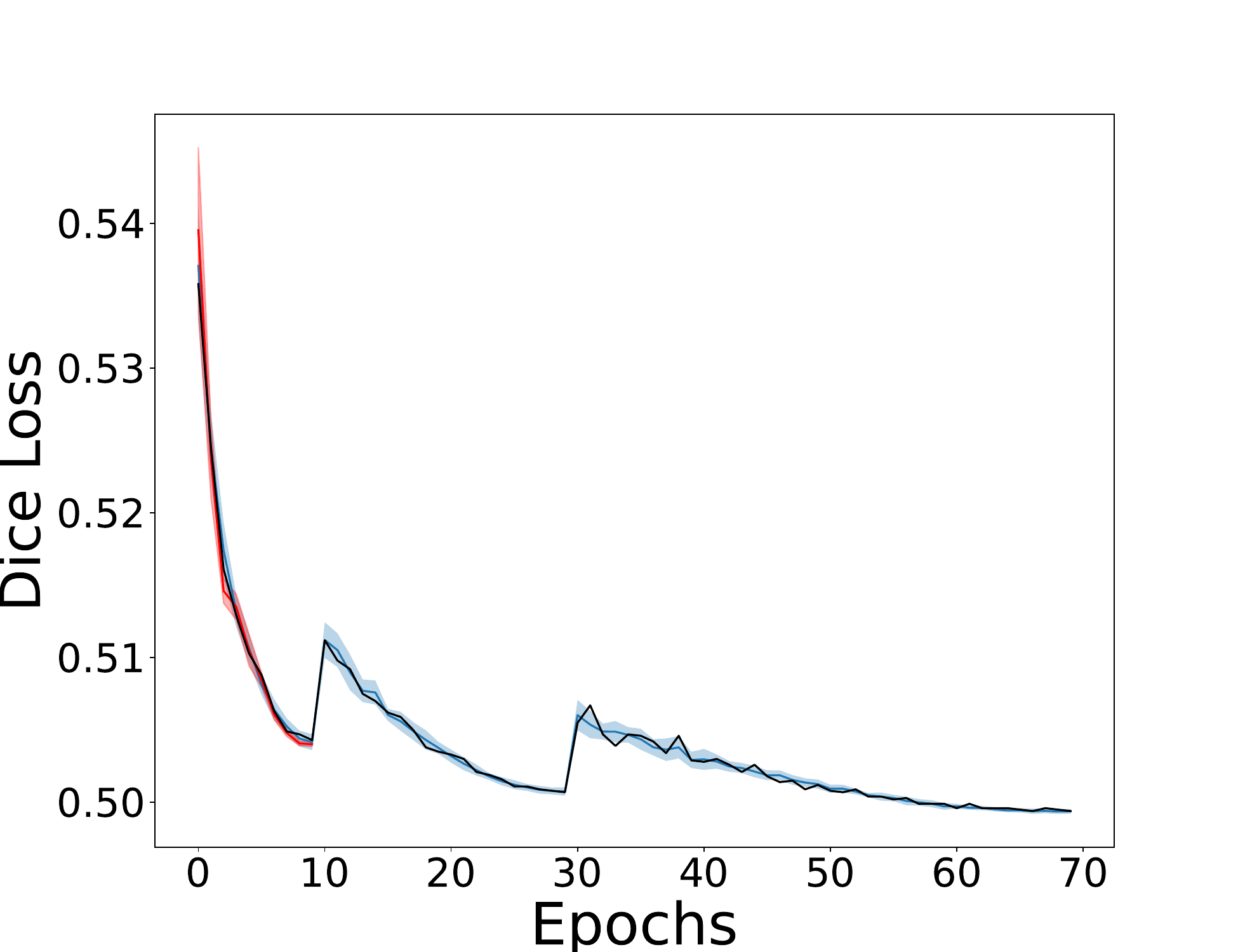}
        \caption{Dice Loss}
        \label{fig:loss_dice}
    \end{subfigure}
    \hspace{-1em}
    \begin{subfigure}[t]{0.33\textwidth}
        \centering
        \includegraphics[width=\linewidth]{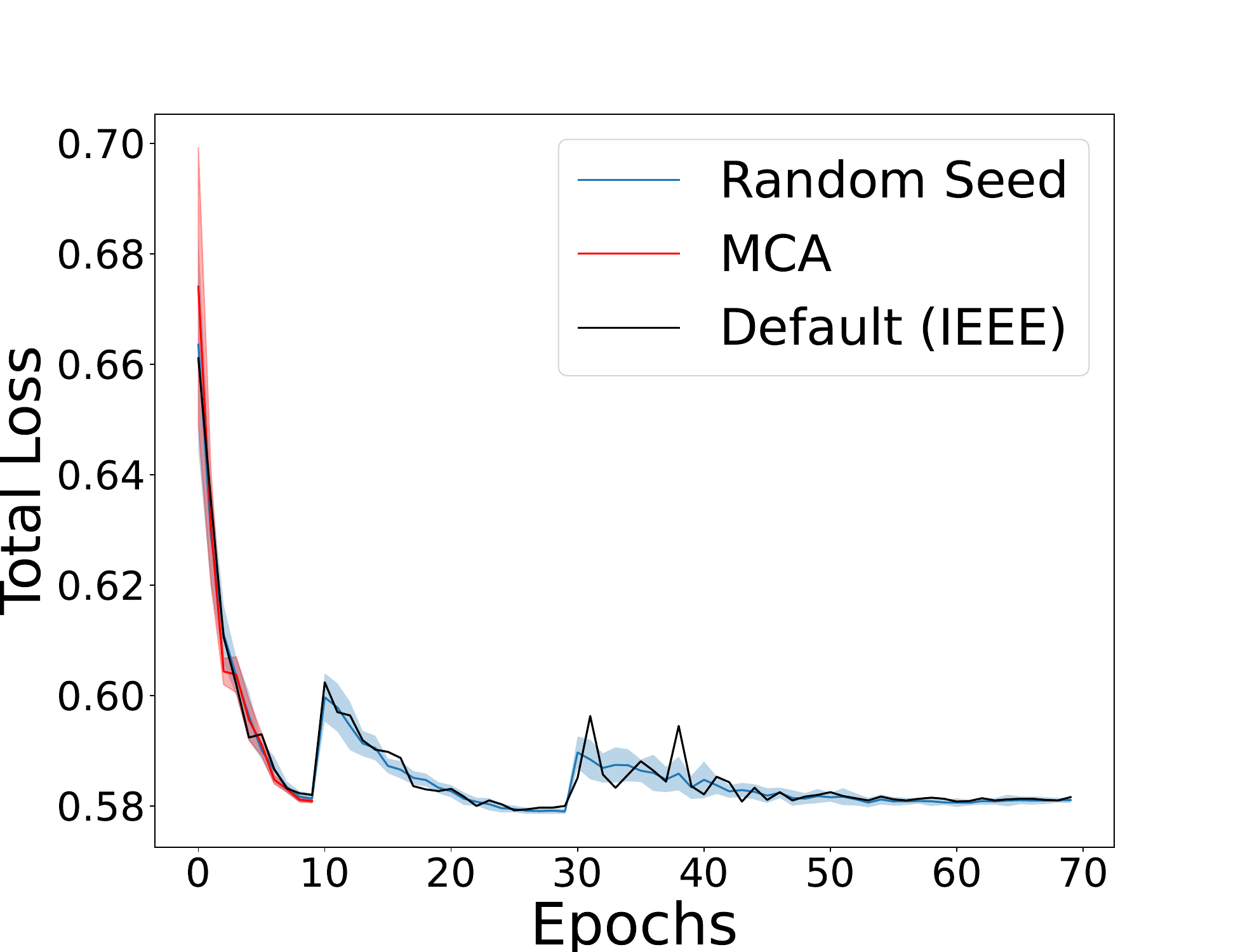}
        \caption{Weighted sum of all losses}
        \label{fig:loss_total}
    \end{subfigure}

    \caption{Validation loss variability for the FastSurfer coronal model. Axial and sagittal models show similar trends. MCA results are evaluated at epoch 10 due to computational overhead; for consistency, Figure~\ref{fig:random_comparison_min_dice} reports all models at this checkpoint. Dice scores at epoch 10 differ by only 4\% from fully trained models, and variability patterns were confirmed to match, indicating that results are representative of the converged solution space.
    }
    \label{fig:fastsurfer_sagittal_loss}
\end{figure}

Figure~\ref{fig:random_comparison_min_dice} compares FastSurfer training
uncertainty (MCA and random seeds) with FreeSurfer MCA uncertainty on the CoRR
subset (left hemisphere shown for clarity; similar trends hold on the right). 
Both FastSurfer MCA and FreeSurfer MCA are conducted on CPU, ensuring hardware-matched conditions for this comparison. Earlier work showed that FastSurfer was substantially less variable than
FreeSurfer at inference for subcortical regions~\citep{pepe2023numerical}. However, during training,
subcortical variability was largely comparable across
methods: for subcortical regions across hemispheres, FastSurfer was significantly less variable
in 3/33 regions, while 12/33 showed higher variability. In contrast, all cortical
regions exhibited significantly greater variability with FastSurfer. 
Overall, this first quantification of numerical uncertainty during training reveals that FastSurfer is more susceptible to variability than its non-DL baseline, suggesting that it is subject to the same reproducibility concerns  identified for FreeSurfer~\citep{chatelain2026practical}. 




Regional variability rankings were strongly correlated between MCA and random seed perturbations (Spearman $\rho = 0.877$, $p = 2.29e-31$). As FastSurfer seed experiments are conducted on GPU, this correspondence spans hardware environments; GPU-specific non-determinism may contribute independently to seed variability, but the strong structural agreement suggests it does not substantially alter the regional pattern of instability. MCA exhibits slightly higher absolute variability overall, as expected: it instruments every floating-point operation and represents a theoretical upper bound, while seeds capture variability through initialization and stochastic dynamics alone.

\subsection{Models Within Ensembles Had Comparable Performance}


We assessed whether training perturbations affected performance relative to a
fixed IEEE baseline (no MCA, fixed seed, single initialization). As shown in
Figure~\ref{fig:fastsurfer_sagittal_loss}, loss trajectories from numerical and
random-seed runs closely tracked the baseline, suggesting that the variability
previously observed could be leveraged through ensembling. Variability increased
primarily after cosine warm restarts ($T_0=10$, $T_{mult}=30$), where learning
rate resets amplified small differences between runs. Despite these
fluctuations, all perturbation methods achieved performance comparable to the
IEEE baseline.
As a sanity check, we replicated the experiment on MNIST under MCA and random seed perturbations (reported on Github), observing consistent patterns as with FastSurfer. 

\subsection{Numerical Ensembling Serves as Data Augmentation Strategy}


\begin{figure}[h]
    \centering
    \includegraphics[width=\linewidth]{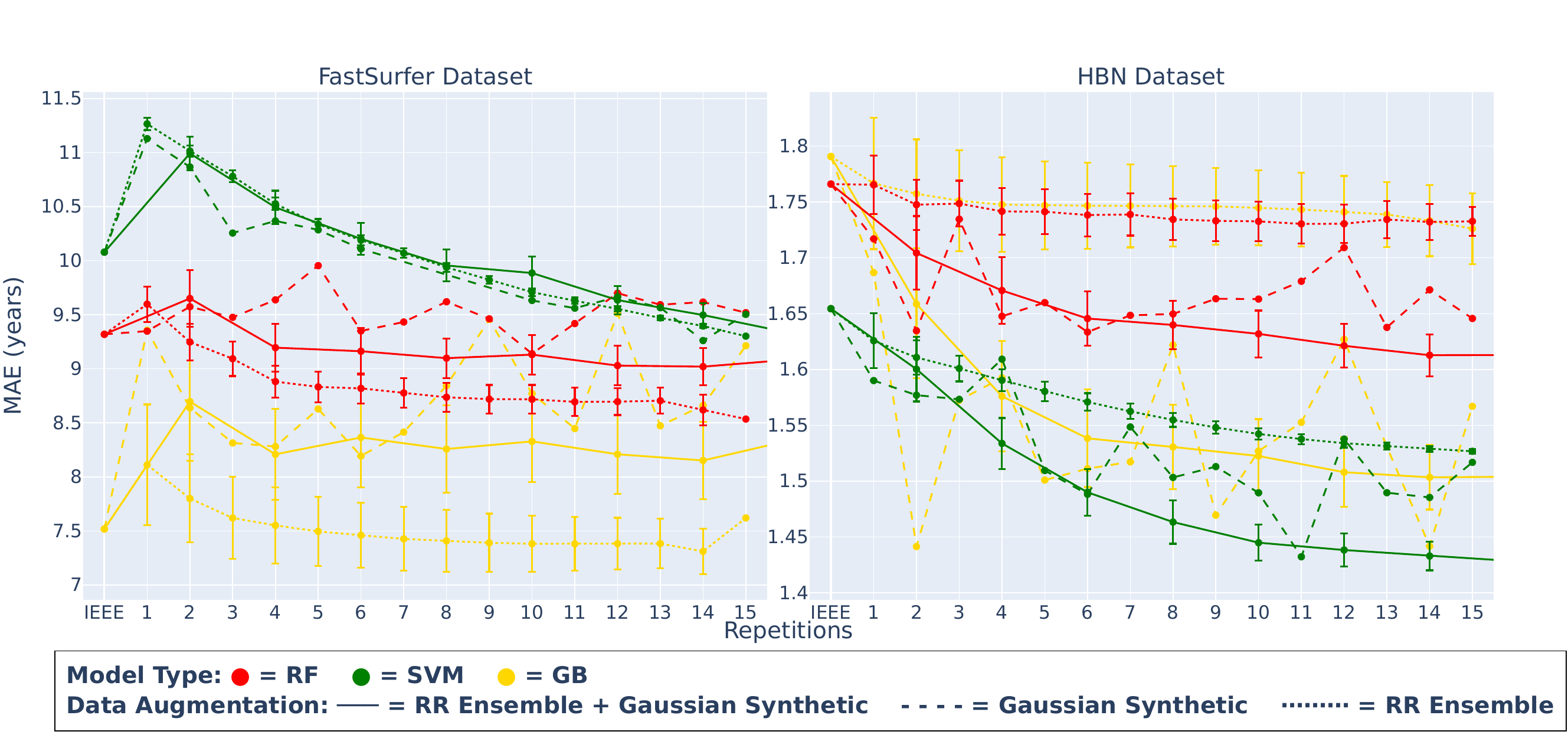}
    \caption{Comparison of brain age regression performance for Random Forest (RF), Support Vector Machine (SVM), and Gradient Boosting (GB) models across the FastSurfer and HBN datasets, highlighting the impact of different data augmentation strategies on mean absolute error (MAE). Chance-level MAE is 20.61 years for FastSurfer and 3.15 years for HBN. 
    }
    \label{fig:data_aug}
\end{figure}

    
    
    



We evaluated data augmentation strategies on both FastSurfer and HBN datasets using RF, SVM, and GB models (Figure~\ref{fig:data_aug}).

\paragraph{FastSurfer Training \& Test Datasets.} Numerical ensembling through repeated segmentations consistently improved performance. MAE decreased significantly with increasing repetitions for SVM ($R=-0.976$, $p<10^{-9}$), RF ($R=-0.849$, $p<10^{-4}$), and GB ($R=-0.671$, $p<10^{-2}$), indicating strong monotonic improvement.
Synthetic augmentation produced weak, non-significant correlations across all models (RF: $R=-0.419$, $p=0.106$; SVM: $R=-0.499$, $p=0.083$); GB: $R=-0.403$, $p=0.122$), while simply adding Gaussian noise yielded no improvement (MAE $\approx$ 24.44 years). Combining approaches showed significant improvement for RF ($R=-0.749$, $p=0.032$) and SVM ($R=-0.974$, $p<10^{-4}$), though GB remained weaker ($R=-0.553$, $p=0.155$).


\paragraph{HBN Dataset.} Ensembling again significantly reduced MAE for RF ($R=-0.832$, $p<10^{-3}$), SVM ($R=-0.973$, $p<10^{-8}$), and GB ($R=-0.908$, $p<10^{-5}$). 
Synthetic augmentation performed comparably to ensembling, but still with weak, non-significant correlations for RF ($R=-0.423$, $p=0.102$), SVM ($R=-0.437$, $p=0.091$), and GB ($R=-0.426$, $p=0.100$), suggesting that HBN's narrower age range reduces the performance gap between augmentation strategies.
Combined strategies showed significant improvement (RF: $R=-0.935$, $p<10^{-4}$; SVM: $R=-0.908$, $p<10^{-3}$; GB: $R=-0.856$, $p<10^{-3}$). GB ensemble performance was more erratic across datasets than RF or SVM, consistent with GB's greater sensitivity to feature noise on a dataset with limited outcome variance.


Notably, absolute MAE differed substantially between datasets: FastSurfer MAE did not fall below 6 years, whereas HBN consistently achieved values below 2 years. This reflects dataset characteristics rather than methodological limitations. FastSurfer spans ages 18-96 years (chance MAE: 20.61), while HBN spans 5-21.9 years (chance MAE: 3.15). In both cases, performance was well below chance.
Our objective was not to achieve state-of-the-art brain age accuracy (typically reported as 2–5 years MAE on cohorts spanning broad adult age ranges; e.g., 20–86 or 50–95 years~\citep{leonardsen2022deep,kumari2024review}), but rather to evaluate numerical ensembling as a principled augmentation mechanism.

Together, these results demonstrate that numerical ensembling provides a strong augmentation strategy. Because each ensemble member is a fully trained FastSurfer model segmenting a real brain MRI, the augmented features reflect authentic neuroanatomical variation. In contrast, synthetic methods must learn distributions from limited data, often producing samples that deviate from realistic anatomy. By remaining anatomically realistic by design, numerical ensembles  outperform synthetic baselines across both datasets.

\section{Conclusion}
We investigated numerical variability in CNN training for neuroimaging, and found that FastSurfer training is at least as sensitive, and in cortical regions more sensitive, to numerical variability than the traditional FreeSurfer pipeline it aims to replace. This finding carries direct implications for downstream analyses: the work in~\citep{chatelain2026practical} recently demonstrated that numerical instability of comparable magnitude in FreeSurfer processing could change clinical conclusions from positive to negative or conversely. Given that FastSurfer training variability meets or exceeds that observed in FreeSurfer, similar risks likely extend to studies relying on FastSurfer-derived measures. 
This motivates the need to systematically investigate uncertainty propagation from DL-based neuroimaging pipelines to downstream analyses, through frameworks such as the Numerical-Population Variability Ratio introduced by~\citep{chatelain2026practical}.

We also found that random seed perturbations produce variability correlated with numerical variability, validating seeds as a practical and well-grounded tool for characterizing this instability. In addition, we show this variability can be harnessed constructively: numerical ensembling via random seeds produces anatomically realistic, diverse representations that improve brain age regression across different datasets, outperforming synthetic baselines without additional data collection.
Together, these results establish that numerical uncertainty in DL training (1) is as substantive as that in classical methods and therefore should be considered when conducting downstream analyses, and (2) can be leveraged for data augmentation.

\bibliography{bib}

\newpage

\appendix

\section{Experimental Reproducibility}
\label{sec:reproducibility}
The code to reproduce this work is located at \url{https://github.com/InesGP/cnn_training_variability}.

\subsection{Computational Infrastructure}
The experiments were conducted on several clusters. Analysis was conducted on a server
equipped with 8 compute nodes with 32 cores Intel(R) Xeon(R) Gold 6130 CPU @
2.10GHz 22MB cache L3. Model training was conducted on 
\begin{itemize}
    \item the Nibi cluster hosted and operated by SHARCNET at University of Waterloo. 
of 134,400 CPU cores and 288 H100 NVIDIA GPUs
    \item the Speed cluster managed by the Concordia Gina Cody School of Engineering and Computer Science
equipped with 24 32-core compute nodes, each with 512 GB of memory, 12 NVIDIA Tesla P6 GPUs, with 16 GB of GPU memory, 
7 servers with 256 CPU cores each 4x A100 80 GB GPUs, partitioned into 4x 20GB MIGs each and 1 AMD FirePro S7150 GPU, with 8 GB of memory
    \item the Virya cluster managed by Concordia University
equipped with 16 x V100 32GB Nvidia GPUs and 8 x A100 40GB Nvidia GPUs (presented as 16 x 20GB MIGs)
\end{itemize}

\subsection{FastSurfer Training Methodology}
\label{sec:fastsurfer_train_steps}

\subsubsection{FreeSurfer} 
FreeSurfer recon-all implements whole-brain segmentation~\citep{fischl2002whole}, through a maximum a-posteriori estimation of the segmentation based on the non-linear registration of the subject image to an atlas by computing a maximum a posteriori (MAP) estimate of the segmentation $W$, given an input MRI scan $I$ and the transform that registers the image into atlas space $L$, where the goal is to maximize $p(W | I,L)$;
\[p(W | I, L) \propto p(I | W,L)~p(W)
\]
$p(I | W,L)$ contains the probability for a given segmentation based off global spatial information and anatomical class statistics about image intensities obtained from the atlas space. $p(W)$ approximates the local spatial relationship between neighboring anatomical classes using anisotropic nonstationary Markov Random Fields.
FreeSurfer is used to preprocess the data given to the FastSurfer model as well as to serve as a reference approach to the model, when both are evaluated on an unseen dataset (CoRR~\citep{zuo2014open}).
We use FreeSurfer v7.3.2, in order to match the version used by the FastSurfer authors during the original model training. 
To measure the variability in FreeSurfer whole brain segmentation, we applied MCA to FreeSurfer using ``fuzzy libmath"~\citep{salari2021accurate}, a version of the GNU mathematical library instrumented with the Verificarlo~\citep{denis2015verificarlo} compiler.

\subsubsection{FastSurfer}
FastSurfer~\citep{henschel2022fastsurfervinn} is a CNN model that performs
whole-brain segmentation, cortical surface reconstruction, fast spherical
mapping, and cortical thickness analysis from anatomical MRIs. The FastSurfer CNN is composed of three 2D fully convolutional
neural networks---each associated with a different 2D slice orientation---that
each have the same encoder/decoder U-net architecture with skip connections,
unpooling layers and dense connections as QuickNAT. A diagram of the model's
architecture is available in the Figure~\ref{fig:fastsurfer}. 
We focus exclusively on the task of whole-brain segmentation, defined as voxel-wise anatomical labeling of brain regions. This segmentation step is entirely performed by the CNN, without surface reconstruction, and uses the pre-trained FastSurfer model (v2.4.0) available on GitHub~\citep{fastsurfer-github}. FastSurfer has demonstrated high accuracy, strong generalization to unseen datasets, and high test-retest reliability.
This model serves as an ideal benchmark for studying numerical variability in high-dimensional medical imaging tasks due to its clinical relevance and architectural complexity.

We closely followed the methodology described in the original FastSurfer paper~\citep{henschel2022fastsurfervinn} and Github repository~\citep{fastsurfer-github}. We obtained access to all training and validation datasets used by the authors (HCP~\citep{VanEssen2013HCP}, ABIDE-I, ABIDE-II, ADNI~\citep{Mueller2005ADNI}, IXI, LA5C, MIRIAD~\citep{Malone2012MIRIAD}, OASIS1~\citep{Marcus2007OASIS1}, OASIS2~\citep{Marcus2009OASIS2}), with the exception of the Rhineland Study dataset, due to data access restrictions, which we replaced with the MICA dataset, which has similar submillimeter resolution. All images were pre-processed using the FreeSurfer HiRes stream with the Desikan-Killiany-Tourville (DKT) atlas. The resulting segmentations were quality-controlled, and cases with failed segmentations were excluded.
For model training, we used the configuration files provided for each anatomical plane, while introducing a perturbation effect across runs. In particular, to investigate the effect of weight initialization, we reinitialized the weights and biases of each layer according to the specific initialization scheme, while keeping all other hyperparameters fixed. Each model was trained for 70 epochs, and the best-performing checkpoint was kept for evaluation. Training typically takes about a day to a day and a half on GPU for the default and random seed models, while weight-initialization experiments can take up to 2.5 days. MCA experiments require over a month, as they must be run on CPU and involve a large computational overhead to perturb floating-point operations.

When comparing the FastSurfer model to the FreeSurfer baseline, we tested both methods on a subset of the Consortium for Reliability and Reproducibility (CoRR) dataset~\citep{zuo2014open}, as it aims to evaluate test-retest reliability and reproducibility and has previously been used for evaluating variability during FastSurfer inference~\citep{pepe2023numerical}

\begin{figure}[H]
  \centering
  \includegraphics[width=\linewidth]{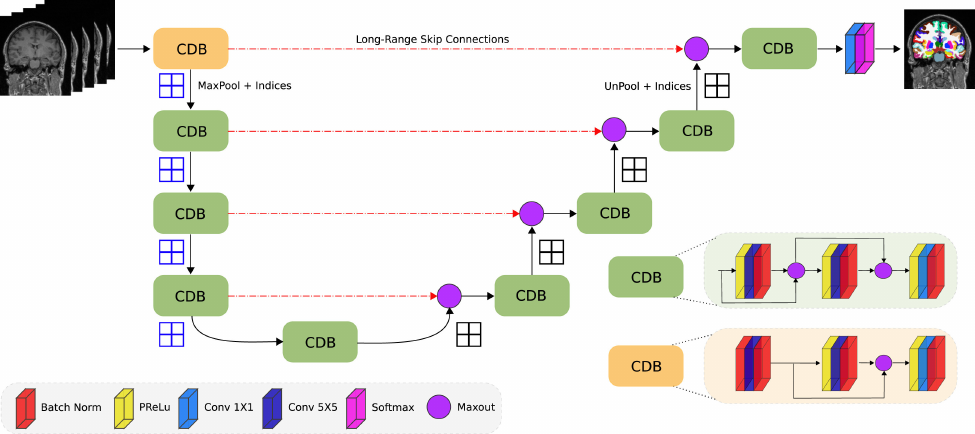}
  \caption{Illustration of FastSurfer's architecture. The CNN consists of four competitive dense blocks (CDB) in the encoder and decoder part, separated by a bottleneck layer. Figure reproduced from~\citep{henschel2020fastsurfer}.}
  \label{fig:fastsurfer}
\end{figure}

\subsection{Data Use Acknowledgment}
\label{sec:data_acks}
Data used in the preparation of this article were obtained from the Alzheimer's Disease Neuroimaging Initiative (ADNI) database (adni.loni.usc.edu). The ADNI was launched in 2003 as a public-private partnership, led by Principal Investigator Michael W. Weiner, MD. The original goal of ADNI was to test whether serial magnetic resonance imaging (MRI), positron emission tomography (PET), other biological markers, and clinical and neuropsychological assessment can be combined to measure the progression of mild cognitive impairment (MCI) and early Alzheimer's disease (AD). For up-to-date information, see adni.loni.usc.edu. 
Data collection and sharing for the Alzheimer's Disease Neuroimaging Initiative (ADNI) is funded by the National Institute on Aging (National Institutes of Health Grant U19AG024904). The grantee organization is the Northern California Institute for Research and Education. In the past, ADNI has also received funding from the National Institute of Biomedical Imaging and Bioengineering, the Canadian Institutes of Health Research, and private sector contributions through the Foundation for the National Institutes of Health (FNIH) including generous contributions from the following: AbbVie, Alzheimer’s Association; Alzheimer’s Drug Discovery Foundation; Araclon Biotech; BioClinica, Inc.; Biogen; BristolMyers Squibb Company; CereSpir, Inc.; Cogstate; Eisai Inc.; Elan Pharmaceuticals, Inc.; Eli Lilly and Company; EuroImmun; F. Hoffmann-La Roche Ltd and its affiliated company Genentech, Inc.; Fujirebio; GE Healthcare; IXICO Ltd.; Janssen Alzheimer Immunotherapy Research \& Development, LLC.; Johnson \& Johnson Pharmaceutical Research \& Development LLC.; Lumosity; Lundbeck; Merck \& Co., Inc.; Meso Scale Diagnostics, LLC.; NeuroRx Research; Neurotrack Technologies; Novartis Pharmaceuticals Corporation; Pfizer Inc.; Piramal Imaging; Servier; Takeda Pharmaceutical Company; and Transition Therapeutics.

Data were provided [in part] by the Human Connectome Project, WU-Minn Consortium (Principal Investigators: David Van Essen and Kamil Ugurbil; 1U54MH091657) funded by the 16 NIH Institutes and Centers that support the NIH Blueprint for Neuroscience Research; and by the McDonnell Center for Systems Neuroscience at Washington University.

LA5C data was obtained from the OpenfMRI database. Its accession number is ds000030
Data were provided [in part] by OASIS-1: Cross-Sectional: Principal Investigators: D. Marcus, R, Buckner, J, Csernansky J. Morris; P50 AG05681, P01 AG03991, P01 AG026276, R01 AG021910, P20 MH071616, U24 RR021382, and OASIS-2: Longitudinal: Principal Investigators: D. Marcus, R, Buckner, J. Csernansky, J. Morris; P50 AG05681, P01 AG03991, P01 AG026276, R01 AG021910, P20 MH071616, U24 RR021382

Data used in the preparation of this article were obtained from the MIRIAD database. The MIRIAD investigators did not participate in analysis or writing of this report. The MIRIAD dataset is made available through the support of the UK Alzheimer's Society (Grant RF116). The original data collection was funded through an unrestricted educational grant from GlaxoSmithKline (Grant 6GKC).

ABIDE data was gathered from various dataset including the California Institute of Technologie funded by the Simons Foundation (SFARI-07-01 to R.A.) and the National Institute of Mental Health (R01 MH080721 to R.A.),  the Carnegie Mellon University funded by NICHD/NIDCD PO1/U19 to M. B. (PI: Nancy Minshew), which is part of the NICHD/NIDCD Collaborative Programs for Excellence in Autism, and the Simons Foundation to M. B. (PI: David Heeger),  the University of Leuven funded by the Flanders Fund for Scientific Research (FWO projects G.0758.10), the FWO postdoctoral Research fellowship grant (KA), the Grant P6/29 from the Interuniversity Attraction Poles program of the Belgian federal government, and the Leuven Autism Researc Consortium (LAuRes), funded by the Research Council of the University of Leuven (IDO/08/013), the Ludwig Maximilians University Munich, the NYU Langone Medical Center funded by NIH (K23MH087770; R21MH084126; R01MH081218; R01HD065282), Autism Speaks, The Stavros Niarchos Foundation, The Leon Levy Foundation, and an endowment provided by Phyllis Green and Randolph Cowen, the Olin, Institute of Living at Hartford Hospital, funded by Autism Speaks (to M.A.), and Hartford Hospital (to M.A) the Social Brain Lab BCN NIC UMC Groningen and Netherlands Institute for Neuroscience funded by Nederlandse Organisatie voor Wetenschappelijk Onderzoek (NWO) Cognitive Pilot Project (051.07.003), NWO VIDI (452-04-305), the NWO Open Competition (400-08-089), the European Commission, Marie Curie Excellence Grant (MEXT-CT-2005-023253), the National Initiative for Brain and Cognition NIHC HCMI Functional Markers (056-13-014), the National Initiative for Brain and Cognition NIHC HCMI Simulation and TOM (056-13-017), and the Netherlands Brain Foundation (KS 2010(1)-29) the Trinity Centre for Health Sciences funded by The Meath Foundation, Adelaide and Meath Hospital, incorporating the National Children's Hospital (AMNCH), Tallaght, and travel fellowship by the Kyulan Family Foundation. the University of Michigan funded by Autism Speaks (CSM), NIH U19 HD035482 and NIH MH066496 (CL), Autism Speaks Pre-doctoral Fellowship 4773 (JLW), the Michigan Institute for Clinical and Health Research (MICHR) Pre-doctoral Fellowship UL1RR024986 (JLW), and NIH R21 MH079871 (SP) the Indiana University funded by National Institutes of Health Grant K99MH094409/R00MH094409 (to D.P.K.) the Erasmus University Medical Center Rotterdam funded by the Simons Foundation Autism Research Initiative (SFARI -307280) and a Dutch ZonMw TOP grant number 91211021 to Tonya White. MRI data acquisition was sponsored in part by the European Community's 7th Framework Programme (FP7/2008-2013 212652). Supercomputing computations were supported by the NWO Physical Sciences Division (Exacte Wetenschappen) and SURFsara (Lisa compute cluster, www.surfsara.nl). The Generation R Study is conducted by the Erasmus Medical Center in close collaboration with the School of Law and Faculty o Social Sciences of the Erasmus University Rotterdam, the Municipal Health Service Rotterdam area, Rotterdam, the Rotterdam Homecare Foundation, Rotterdam and the Stichting Trombosedienst \& Artsenlaboratorium Rijnmond (STAR-MDC), Rotterdam.

The IXI data was obtained from their website at https://brain-development.org/ixi-dataset/


\section{MNIST Benchmark}
\label{sec:mnist}


\begin{figure}[htbp]
    \centering
    \includegraphics[width=\linewidth]{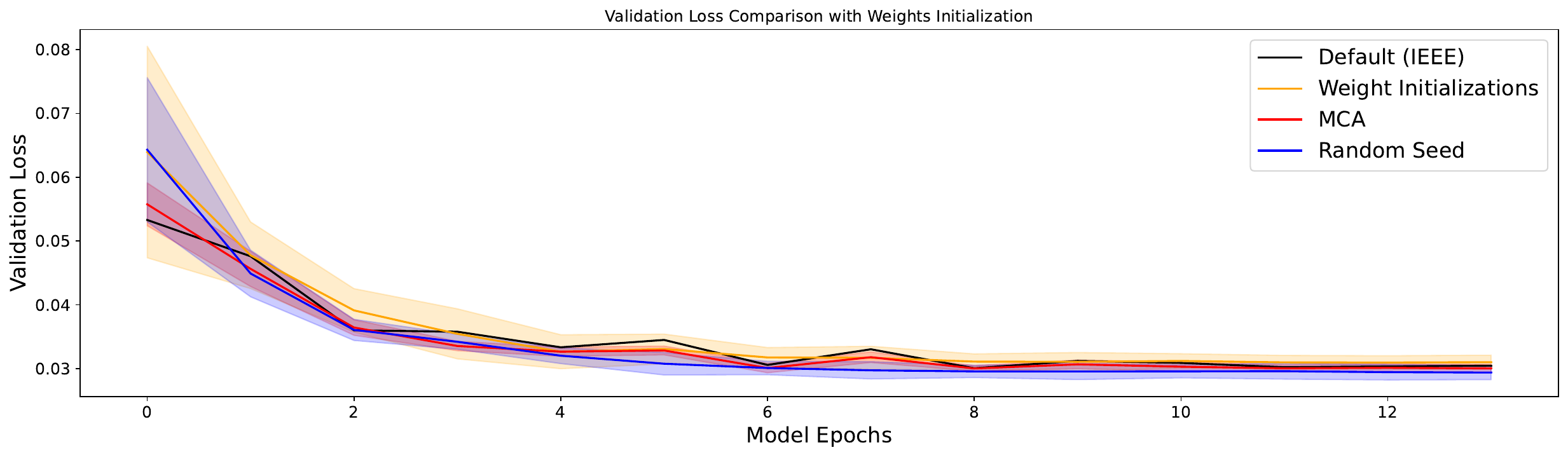}
    \caption{Comparison of sources of variability across validation loss for MNIST model}
    \label{fig:mnist_val_loss}
\end{figure}

\begin{figure}[htbp]
    \centering
    \includegraphics[width=\linewidth]{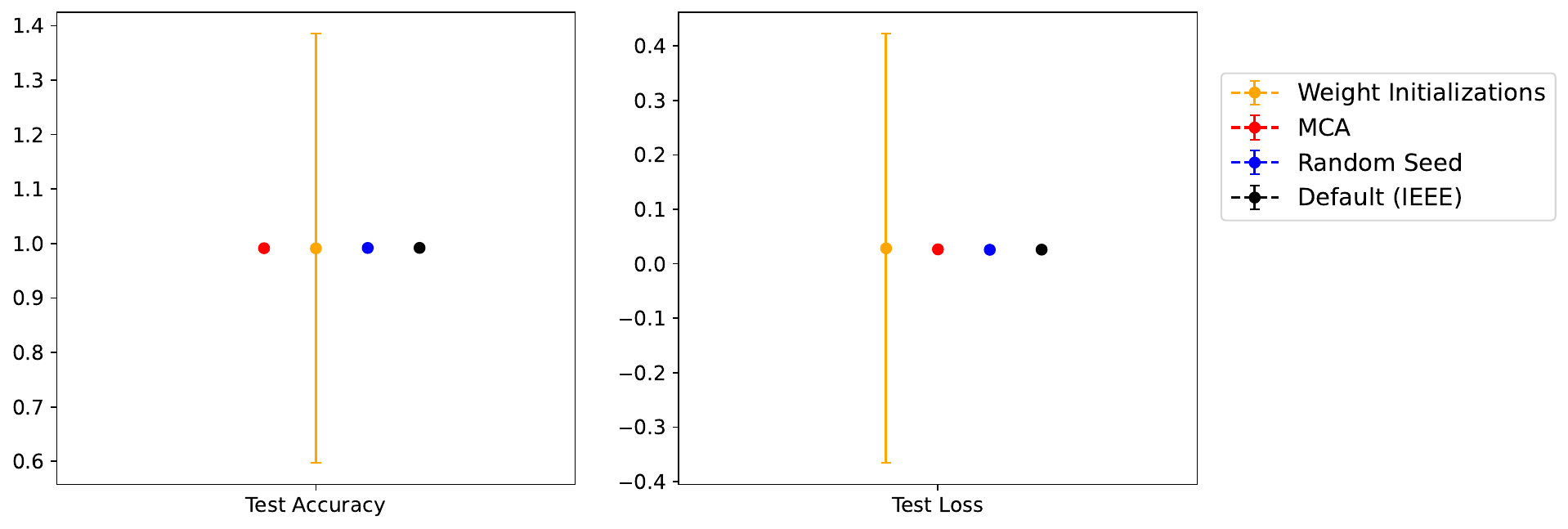}
    \caption{Comparison of sources of variability across test accuracy and loss for MNIST model}
    \label{fig:mnist_test}
\end{figure}

To generalize beyond FastSurfer, we repeated these experiments on a CNN trained on MNIST. Here, weight initialization again dominated variability, with differences on the order of $10^{-1}$ in validation loss (Figure~\ref{fig:mnist_val_loss}). MCA and random seed perturbations remained comparatively minor ($\sim10^{-3}$), though random seeds produced larger deviations ($\sim10^{-2}$) early in training—consistent with their influence on initialization and dropout. Figure~\ref{fig:mnist_test} highlights the performance impact: while weight initialization reduced test accuracy to as low as 94\%, both MCA and random seed variability maintained near-baseline accuracy ($\sim99\%$).

\end{document}